\documentclass[12pt]{article}

\usepackage[T2A]{fontenc}
\usepackage[utf8]{inputenc}

\usepackage{amsfonts}
\usepackage{amssymb}
\usepackage{amsmath}
\usepackage{amsthm}

\def \le {\leqslant}
\def \ge {\geqslant}

\begin{document}

\begin{Large}
\centerline{\bf On Kleinbock's  Diophantine result }
\end{Large}
 \vskip+0.5cm \centerline{by {\bf Nicky Moshchevitin}.\footnote{supported by the grant RFBR № 09-01-00371
 }  }

 \vskip+0.5cm
\begin{small}
{\bf Abstract.}\,
We give an elementary proof of a  
metrical Diophantine 
result by D. Kleinbock realted to  badly approximable vectors in affine subspaces.
\end{small}

\vskip+0.5cm
 
{\bf 1.\, Definitions and notation.}
\,\,\,
Let
$\mathbb{R}^d$
be a Euclidean space with the coordinates $(x_1,...,x_d)$,
let $\mathbb{R}^{d+1}$
be a Euclidean space with the coordinates $(x_0,x_1,...,x_d)$.
For ${\bf x} \in \mathbb{R}^d$ or
${\bf x} \in \mathbb{R}^{d+1}$ we denote by
$|{\bf x}|$ its sup-norm:
$$
|{\bf x}|=
\max_{1\le j\le d} |x_j|\,\,\,\,
\text{or}\,\,\,\,
|{\bf x}|=
\max_{0\le j\le d} |x_j|
.
$$
Consider an affine subspace $A\subset \mathbb{R}^d$ such that ${\bf  0}\not\in A$ and define
the affine subspace
${\cal A}\subset \mathbb{R}^{d+1}$ in the following way:
$$
{\cal A} =\{ {\bf x} = (1,x_1,...,x_d):\,\, (x_1,...,x_d) \in A\}.
$$
Let $B$ be another affine subspece such that
$
B\subset A$. Define
$$
{\cal B} =\{ {\bf x} = (1,x_1,...,x_d):\,\, (x_1,...,x_d) \in B\},\,\,\,\,\,
{\cal B} \subset {\cal A}.
$$
Put
$$
a ={\rm dim } \,A= {\rm dim } \,{\cal A},\,\,\,\, b = {\rm dim }\, B= {\rm dim } \,{\cal B}.
$$
We define {\it linear} subspaces
$$
\hbox{\got A} = {\rm span} \,{\cal A},\,\,\,\,\
\hbox{\got B} = {\rm span }\,{\cal B}
$$
as the smallest linear subspaces in $\mathbb{R}^{d+1}$ containing
${\cal A}$ or ${\cal B}$ respectively.
So
$$
{\rm dim } \, \hbox{\got A}  =a+1,\,\,\,\,{\rm dim }\, \hbox{\got B}= b+1. 
$$

Let $\psi (T), \,\,T\ge 1$ be a positive valued function decreasing to zero as $ T \to +\infty$.
We define an affine  subspace $B$ to be {\it $\psi$-badly approximable} if
\begin{equation}
\inf_{{\bf x} \in \mathbb{Z}^{d+1}\setminus \{{\bf 0}\}}\,\left(
\frac{1}{\psi (|{\bf x}|)}\,\inf_{ {\bf y} \in \hbox{\got B}}|{\bf x}-{\bf y}|\right)>0.
\label{def}
\end{equation}
This definition is very convenient for our exposition. 

Here 
we would like to note that in the case when the affine subspace $B$ has zero dimension
(and hence $B=\{ {\bf w}\}$ consists of just one nonzero vector
${\bf w} = (w_1,....,w_d)  \in \mathbb{R}^d$)
the definition (\ref{def}) gives
\begin{equation}
\inf_{{\bf x} \in \mathbb{Z}^{d+1}\setminus \{{\bf 0}\}}\,\left(
\frac{1}{\psi (|{\bf x}|)}\,\inf_{t\in \mathbb{R}}|{\bf x}-t{\bf w}^*|\right)>0
\label{def1}
\end{equation}
with ${\bf w}^* =(1,w_1,...,w_d)$.
Under mild restrictions   on the function $\psi$
 the inequality (\ref{def1}) is equivalent to the condition
that there exists a positive constant $\gamma = \gamma ({\bf w})$ such that
\begin{equation}
\max_{1\le j \le d} ||w_j  q || \ge \gamma \cdot \psi (|q|),\,\,\,\,
\forall \, q \in \mathbb{Z} \setminus \{ 0\}
\label{bad}
\end{equation}
(here $||\cdot ||$ denotes the distance to the nearest integer).
We can consider the inequality 
\begin{equation}\label{restric}
 \inf_{T\ge 1} \frac{\psi (\kappa T)}{\psi (T)} > 0, \,\,\,\,\,\, \forall \, \kappa \ge 1
\end{equation}
as an example of such a condition on the function $\psi$.

The condition (\ref{bad}) is a usual condition of $\psi$-badly simultaneously approximability 
of the vector ${\bf w}$.
This explains our definition (\ref{def}).

In the sequel we consider two positive decreasing 
to zero 
functions $\psi (T), \varphi (T)$ under the condition
$
\varphi (T)\le \psi (T)
$.
We  suppese that $\varphi (T)$ satisfies the condition
\begin{equation}\label{rend}
 \sup_{T\ge 1} \frac{\varphi (\kappa T)}{\varphi (T)} 
\le +\infty , \,\,\,\,\,\,\,\ \forall \, \kappa \ge 1
.\end{equation}
For these functions we  define quantities
$$\mu_T =
\left(\frac{
T}{\psi (T)}\right)^{a-b}
,\,\,\,\,\,
\lambda_T =
\left(
\frac{\varphi(T)}{T}
\right)^{a}-\,\,
\left(
\frac{\varphi(T+1)}{T+1}
\right)^{a}
.
$$

{\bf 2. \, The result.}
\,\,\,
Now we are ready to formulate the main result of this paper.

{\bf Theorem 1.}\,\,{\it  
 Suppose that  for all $T\ge 1$ one has
$0<\varphi (T) \le \psi (T) $ and
the series
\begin{equation}
\sum_{T=1}^{+\infty}
\mu_T
\lambda_T
\label{ser}
\end{equation}
converges.
Suppose that $\psi$  satisfies (\ref{restric}) and $\varphi$ satisfies (\ref{rend}).
Suppose that $B\subset A$, ${\bf 0}\not\in A$ are affine subspaces and $ 0\le b={\rm dim}\, B< a={\rm dim}\, A\le d$.
 Let $B$ be a $\psi$-badly approximable affine subspace.

 Then  almost all
(in the sense of Lebesgue measure)
 vectors ${\bf w}\in A$ 
 are $\varphi$-badly approximable vectors.
}

We consider a special case  $b=0$ and $\psi_d (T) = T^{-1/d}$.
Then $\psi_d$-badly approximable vectors ${\bf w}$
are known as simultaneously badly approximable vectors (see \cite{SC}, Chapter 2).
Take
\begin{equation}\label{DF}
\varphi_{u,\Delta} (T) = \psi_d (T) \cdot (\log T)^{-\Delta} =
T^{-\frac{1}{d}} \, (\log T)^{-\Delta}.
\end{equation}
In the case 
$
\Delta
>\frac{1}{a}
$
the series (\ref{ser}) converges. 
So we have the following

{\bf Corollary 1.}\,\,\,{\it Consider  an affine subspace $A\subset \mathbb{R}^d$, ${\bf 0}\not\in A$. Suppose  that there exists a badly approximable
vector ${\bf w} \in A$.  Suppose that $\Delta > \frac{1}{a}$. Then almost all vectors from the subspace $A$ are
$\varphi_{d,\Delta}$-badly approximable vectors,
where $\varphi_{d,\Delta} (T)$ is defined in (\ref{DF}).}

Here we should note that in the case $b =1$ and $ \psi (T) = T^{-1/d}$ such a result was obtained by
Dmitry Kleinbock  (see \cite{K}, Theorem 4.2)
by means of theory of dynamics on homogeneous spaces.
Some generalizations are due to 
Dmitry Kleinbock  \cite{K1}
and
Yuqing Zhang \cite{z}.
We would like
to note that the main subject of the papers \cite{K,K1,z} is dynamical approach to Diophantine 
approximations on submanifolds.
A weaker result was announced  by Mikhail B. Sevryuk \cite{sev}.

Consider the function $\psi_a(T) = T^{-\frac{1}{a}}$.
One can easily see that for a given affine subspace $A\subset \mathbb{R}^d$ 
of dimension ${\rm dim } \,A = a\ge 1$
there exists a $\psi_a$-badly approximable vector ${\bf w}$ such that ${\bf w} \in A$.
 Obiously we have the following  

{\bf Corollary 2.}\,\,\,{\it For  any affine subspace $A\subset \mathbb{R}^d$, ${\bf 0}\not\in A$
almost all vectors
 from $A$ are
$\varphi_{a,\Delta}$-badly approximable vectors,
where $\varphi_{a,\Delta} (T)$ is defined in (\ref{DF}).}

{\bf 3. \, Proof of Theorem 1.}
\,\,\, Take $R\ge 1$.
 In the sequel we do not take care on  the constatns. All constants in the symbols $\ll,\asymp$ may depend on $d$, subspaces $A,B$.
and $R$.

For a set $\hbox{\got C}\subset \mathbb{R}^{d+1}$ and a point ${\bf x}\in \mathbb{R}^{d+1}$ we define the distance
$|{\bf x}|_{\hbox{\got C}}$
from ${\bf x}$ to $\hbox{\got C}$ by
$$
|{\bf x}|_{\hbox{\got C}}=
\inf_{{\bf y} \in \hbox{\got C}}|{\bf x} -{\bf y}|.
$$ 

Consider the set
$$
\Omega_T =\{{\bf z} =(z_0, z_1,...,z_d)\in \mathbb{R}^{d+1}:\, 
0\le z_0\le T,\,\,\,
\max_{1\le j \le d}|z_j|\le RT,
\,\,\, |{\bf z}|_{\hbox{\got B}}\le \gamma \cdot \psi (RT)
\}.
$$
Here we suppose that
$\gamma >0$ is strictly less than the infimum in (\ref{def}).
As $B$ is  a $\psi$-badly approximable subspace we see that
$$
\Omega_T \cap \mathbb{Z}^{d+1} = \{{\bf 0}\}
$$
(we  use the definition (\ref{def})).
Now we observe that any translation of the $1/2$-dilatated set
\begin{equation}
\frac{1}{2}\cdot \Omega_T + {\bf c},\,\,\,\, {\bf c} \in \mathbb{R}^{d+1}
\label{half}
\end{equation}
consists of not more than one integer point.
Indeed if  two different integer points $ {\bf x}, {\bf y}$ belong to the same set of the form (\ref{half}) then
${\bf 0}\neq {\bf x} -{\bf y} \in \Omega_T$. This is not possible.

 Consider the set
$$
\Pi_T =
\{
{\bf z}\in \mathbb{R}^{d+1}:\, 0\le z_0\le T,\,\,\,
\max_{1\le j \le d}|z_j|\le RT,\,\,\, |{\bf z}|_{\hbox{\got A}}\le\varphi (RT)
\}.
$$
Let $\nu_T$ be the minimal number of  points ${\bf c}_j, 1\le j \le \nu_T$ such that
$$
\Pi_T \subset \bigcup_{j=1}^{\nu_T}\left(
\frac{1}{2}\cdot \Omega_T + {\bf c}_j \right).
$$ 
As $\varphi (T)\le \psi (T)$ and $\psi$ satisfies (\ref{restric}), we see that
the  set $
\Pi_T $ can be covered by 
$\nu_T\ll
\mu_T
$
different sets of the form (\ref{half}).
 So we  deduce a conclusion about an upper bound for the number of integer points in 
$\Pi_T$:
\begin{equation}
 \# \left(\Pi_T\cap \mathbb{Z}^{d+1}\right) \ll 
 \mu_T.
\label{int}
\end{equation}

For an integer $T\ge 1$  we consider the set 
$$
 \mathbb{Z}_T =\{{\bf z}=(z_0,z_1,...,z_d)\in \mathbb{Z}^{d+1}:\,\, z_0= T,\,\,\,
\max_{1\le j \le d}|z_j|\le RT,\,\,\,\,
 |{\bf z}|_{\hbox{\got A}}\le \varphi (RT)\}
$$
of the cardinality
\begin{equation}
\zeta_T^{(R)} =\# \mathbb{Z}_T,\,\,\,\,\,\,
0\le  \zeta_T^{(R)}\ll T^a.
\label{ze}
\end{equation}
For $\rho >0$ and ${\bf z}\in \mathbb{R}^{d+1}$ put 
$$
\hbox{\got U}_\rho ({\bf z}) =
\{ {\bf y}\in \mathbb{R}^{d+1}:\,\,\, |{\bf z}-{\bf y}|\le \rho\}.
$$
Consider the set
$$
\hbox{\got U}_T =
\bigcup_{{\bf z}} \hbox{\got U}_{\varphi (RT)} ({\bf z}) 
$$
where  the union is taken over all integer points ${\bf z}$ such that
$$z_0= T,\,\,\,\,\,\,\,
\max_{1\le j \le d}|z_j|\le RT,\,\,\,\,\,\,\,
\hbox{\got U}_{\varphi (RT)} ({\bf z})\cap \hbox{\got A}\neq \varnothing.$$
Clearly
$$
\hbox{\got U}_T =
\bigcup_{{\bf z}\in \mathbb{Z}_T} \hbox{\got U}_{\varphi (RT)} ({\bf z}) 
$$
Consider the cone
$$
\hbox{\got G}_T =\{
{\bf x}\in \mathbb{R}^{d+1}:\,\, {\bf x} = t\cdot {\bf y},\,\, t\in \mathbb{R}, {\bf y } \in 
\hbox{\got U}_T\}
$$
and the projection
$$
{\cal U}_T = {\cal A}\cap  \hbox{\got G}_T .
$$
Consider the series
 \begin{equation}
\sum_{T=1}^{+\infty}
{\rm mes}_{a}\,
{\cal U}_T
\label{newser}
\end{equation}
where  ${\rm mes}_{a}$ stands for the $a$-dimensional Lebesgue measure.
Suppose that  the series (\ref{newser}) converges. Then by the Borel-Cantelli lemma
argument 
we see that the set of non-$\varphi$-badly approximable vectors ${\bf w}\in {\cal A}\cap \{ |{\bf z}| \le R\}$
is a set of zero measure. As $R$ is arbitrary we see that
Theorem 1 follows from the 
 convergence of the series (\ref{newser}).

Now we show that the series (\ref{newser}) converges.
Note that
the set
$\hbox{\got U}_T$
is a union of not more than
$
\zeta_T^{(R)}
$
 balls
(in sup-norm)
 of the radius
$\varphi (RT)$.
Hence the set ${\cal U}_T\subset {\cal A}$ can be covered by not more than $\zeta_T^{(R)}$ balls of the radius
$\varphi (RT)/T$.
So in order to prove the convergence of the series (\ref{newser}) one can establish the convergence of the series
\begin{equation}
\sum_{T=1}^{+\infty}
\zeta_T^{(R)}\cdot
\left(\frac{\varphi (T)}{T}\right)^a
\label{newnewser}
\end{equation}
(we take into account (\ref{rend})).

It follows from (\ref{int}) and monotonicity of $\psi$ that
$$
 \sum_{j=1}^T\zeta_j^{(R)}\ll
\sum_{\nu \le \log T} \mu_{T/2^\nu}
=
\sum_{\nu \le \log T} \left(\frac{T/2^\nu}{\psi (T/2^\nu)} \right)^{a-b}\ll
\left(\frac{T}{\psi (T)} \right)^{a-b}
 =\mu_T.
$$
Now we  see that
the convergence of the series (\ref{newnewser}) follows from the convergence of the series (\ref{ser})
 by partial summation as
from (\ref{ze}) we see that
$$
 \zeta_T^{(R)}\cdot \left(\frac{\varphi (T+1)}{T+1}\right)^a
\ll (\varphi (T))^a
\to 0,\,\,\,\,\,
T\to +\infty.
$$
Theorem 1 is proved.

\end{document}